\documentstyle[11pt,leqno]{article}
\newtheorem{theorem}{{\sc Theorem}}
\newcommand{\bt}{\begin{theorem}}
\newcommand{\et}{\end{theorem}}
\newcommand{\newsection}[1]{\setcounter{equation}{0} \setcounter{theorem}{0}
\section{#1}}

\newcommand{\NI}{\noindent}
\newcommand{\bea}{\begin{eqnarray}}
\newcommand{\eea}{\end{eqnarray}}

\def \spec#1 {\mathop{#1}}

\def \b #1 {\bf #1}

\newcommand{\ity}{\infty}
\newcommand{\raro}{\rightarrow}

\newcommand{\vsp}{\vskip 1em}

\newcommand{\be}{\begin{equation}}
\newcommand{\ee}{\end{equation}}
\newcommand{\ben}{\begin{eqnarray*}}
\newcommand{\een}{\end{eqnarray*}}

\begin{document}
\NI{\bf On a characterization of probability distribution based on}\\
\NI{\bf maxima of independent or max-independent random variables}\\
\vsp
\NI{\bf B.L.S. Prakasa Rao\footnote{CONTACT B.L.S. Prakasa Rao, blsprao@gmail.com, CR Rao Advanced Institute of Mathematics, Statistics and Computer Science, Hyderabad 500046, India}}
\vsp
\NI{\bf CR Rao Advanced Institute of Mathematics, Statistics and  Computer Science, Hyderabad, India}
\vsp
\NI{\bf Abstract:}\\ Kotlarski (1978) proved a result on identification of the distributions of independent random variables $X,Y$ and $Z$ from the joint distribution of the bivariate random vector $(U,V)$ where $(U,V)= (\max (X,Z), \max(Y,Z)).$ We extend this result to the case $$(U,V)=(\max(X,aZ_1,bZ_2),\max(Y,cZ_1,dZ_2)$$
where $X,Y,Z_1,Z_2$ are independent or max-independent random variabkes, $Z_1$ and $Z_2$ are identically distributed and $a,b,c,d$ are known positive constants.
\vsp
\NI{\bf Keywords:}\\ Kotlarski's lemma; Identifiability; Characterization; Max-independent; Independent; Maxima.  

\vsp
\NI{\bf MSC 2020: Primary 62E10}\\
\vsp
\newsection{\bf Introduction} 
\vsp
Let $X_0,X_1$ and $X_2$ be independent random variables. Define $Y_1= \max(X_0,X_1)$ and $Y_2= \max(X_0,X_2)$. It is of interest to know whether the joint distribution of $(Y_1,Y_2)$ determines the individual distributions of $X_0,X_1$ and $X_2$ uniquely. It is known that the random variable $Y_1$ alone can not determine the distributions of $X_0$ and $X_1$ uniquely unless $X_0$ and $X_1$ are identically distributed random variables (cf. Prakasa Rao (1992), Section 7.3).  Kotlarski (1978) and Klebanov (1973) obtained characterizations for probability distributions  through maxima or minima of independent random variables. Prakasa Rao (2024) discussed characterizations of probabilty distributions based on maxima or minima of some families of dependent random variables. We now discuss extension of the result in Kotlarski (1978)  leading to characterizations of probability distributions through maxima for some classes of independent or max-independent random variables.
\vsp
\newsection{\bf Identifiability by maxima for independent random variables}
\vsp
The following result is due to Kotlarski (1978).
\vsp
\NI{\bf Theorem 2.1:} (Identifiability by maxima) Suppose $X_0,X_1,X_2$ are independent random variables. Define $Y_1= \max(X_0,X_1)$ and $Y_2= \max(X_0,X_2)$. Then the joint distribution of $(Y_1,Y_2)$ uniquely determines the distributions of the independent random variables $X_0,X_1$ and $X_2$ provided the supports of their distribution functions are the same.
\vsp
For a proof of Theorem 2.1, see Kotlarski (1978) (cf. Prakasa Rao (1992), Theorem 2.2.1, p.24).
\vsp
We will now generalize Theorem 2.1 in analogy with the results of Li and Zheng (2019) for linear functions of independent random variables.
\vsp
\NI{\bf Theorem 2.2:} Let $U$ and $V$ be random variables defined by the relations
$$U=\max(X,aZ_1,bZ_2);\;V=\max(Y,cZ_1,dZ_2)$$
where $X,Y,Z_1$ and $Z_2$ are independent random variables, $Z_1,Z_2$ are identically distributed and $a,b,c,d$ are known positive constants. Further suppose that the distribution functions of $X,Y,Z_1$ have the same support. Then, the joint distribution function of $(U,V)$ uniquely determines the distributions of $X,Y,Z_1$ if $a=b$ or if $a \neq b$ but the distribution function of $Z_1$ is differentiable and the derivative is continuous.
\vsp
\NI{\bf Proof:} Let $G(t_1,t_2)$ be the joint distribution function of $(U,V).$ Suppose that $(N,M, S_1,S_2)$ is another set of independent random variables with $S_1$ and $S_2$ identically distributed with the same support as that of $X$ and such that the joint distribution of the random vector
$$(\max(N,aS_1,bS_2), \max(M,cS_1,dS_2))$$
is the same as that of the random vector $(U,V).$ Let $F_X(.)$ denote the distribution function of the random variable $X.$ It can be seen that, for any$ -\ity<t_1,t_2<\ity,$
\bea
G(t_1,t_2)&=& P (\max(X,aZ_1,bZ_2)\leq t_1,\max((Y,cZ_1,dZ_2)\leq t_2)\\\nonumber
&=& P(X \leq t_1,aZ_1\leq t_1,bZ_2\leq t_2; Y\leq t_2,cZ_1\leq t_2, dZ_2\leq t_2)\\\nonumber
&=& P(X\leq t_1, Y\leq t_2, Z_1\leq \frac{t_1}{a}, Z_1\leq \frac{t_2}{c}, Z_2\leq \frac{t_1}{b}, Z_2\leq \frac{t_2}{d})\\\nonumber
&=& P(X\leq t_1, Y\leq t_2, Z_1\leq \min(\frac{t_1}{a},\frac{t_2}{c}), Z_2\leq \min(\frac{t_1}{b},\frac{t_2}{d}))\\\nonumber
&=& F_X(t_1)F_Y(t_2)F_{Z_1}(\min(\frac{t_1}{a},\frac{t_2}{c}))F_{Z_1}(\min(\frac{t_1}{b},\frac{t_2}{d}))\\\nonumber
\eea
since $X,Y Z_1,Z_2$ are independent random variables and $Z_1,Z_2$ are identically distributed. By a similar argument it follows that
\bea
G(t_1,t_2)&=& P (\max(N,aS_1,bS_2)\leq t_1,\max(M,cS_1,dS_2)\leq t_2)\\\nonumber
&=& F_N(t_1)F_M(t_2)F_{S_1}(\min(\frac{t_1}{a},\frac{t_2}{c}))F_{S_1}(\min(\frac{t_1}{b},\frac{t_2}{d}))\\\nonumber
\eea
Hence
\bea
F_X(t_1)F_Y(t_2)F_{Z_1}(\min(\frac{t_1}{a},\frac{t_2}{c}))F_{Z_1}(\min(\frac{t_1}{b},\frac{t_2}{d}))\\\nonumber
=F_N(t_1)F_M(t_2)F_{S_1}(\min(\frac{t_1}{a},\frac{t_2}{c}))F_{S_1}(\min(\frac{t_1}{b},\frac{t_2}{d})) \\\nonumber 
\eea
for all $-\ity <t_1, t_2 <\ity.$ Let $D\subset R$ be the common support of $X,Y,Z_1$ and $N,M,S_1.$ Then, for all $t\in D,$
define
$$\eta_1(t)= \frac{F_X(t)}{F_M(t)}, \eta_2(t)= \frac{F_Y(t)}{F_N(t)}\;\;\mbox{and}\;\; \eta_3(t)= \frac{F_{Z_1}(t)}{F_{S_1}(t)}.$$
Equation (2.3) implies that
\be
\eta_1(t_1)\eta_2(t_2)\eta_3(\min(\frac{t_1}{a},\frac{t_2}{c}))\eta_3(\min(\frac{t_1}{b},\frac{t_2}{d}))=1,t_1, t_2 \in D. 
\ee
Let $t_1\raro \ity$ in equation (2.4). From the properties of the distribution functions, it follows that 
$$\eta_2(t_2)\eta_3(\frac{t_2}{c})\eta_3(\frac{t_2}{d})=1, t_2\in D.$$
Hence
\be
\eta_2(t_2)=[\eta_3(\frac{t_2}{c})\eta_3(\frac{t_2}{d})]^{-1}, t_2\in D. 
\ee
Let $t_2 \raro \ity$ in equation (2.4). From the properties of the distribution functions again, it follows that 
$$\eta_1(t_1)\eta_3(\frac{t_1}{a})\eta_3(\frac{t_1}{b})=1, t_1\in D$$ 
which implies that
\be
\eta_1(t_1)=[\eta_3(\frac{t_1}{a})\eta_3(\frac{t_1}{b})]^{-1}, t_1\in D. 
\ee
Combining the equations (2.4)-(2.6), it follows that
\be
\eta_3(\min(\frac{t_1}{a},\frac{t_2}{c}))\eta_3(\min(\frac{t_1}{b},\frac{t_2}{d}))= \eta_3(\frac{t_1}{a})\eta_3(\frac{t_1}{b})\eta_3(\frac{t_2}{c})\eta_3(\frac{t_2}{d}), t_1,t_2\in D.
\ee
Fix $t_1 \in D$ and let $t_2 \raro \ity$ or $t_2$ tend to the right boundary point of the set $D$. Then the expression on the left side of equation (2.7) tends to 1 and the expression on the right side tends to 
$$\eta_3(\frac{t_1}{a})\eta_3(\frac{t_1}{b})$$
by the properties of the distribution functions. Hence
\be
\eta_3(\frac{t_1}{a})\eta_3(\frac{t_1}{b})=1, t_1\in D. 
\ee
This in turn implies that
\be
\eta_1(t_1)=1, t_1\in D 
\ee
from (2.6). A similar analysis shows that 
\be
\eta_2(t_2)=1, t_2\in D 
\ee
by fixing $t_2\in D$ and letting $t_1\raro \ity$ or $t_1$ tend to the right boundary point of the set $D$.
In particular, it follows that
$$F_X(t)=F_M(t), t \in D$$
and
$$F_Y(t)=F_N(t), t \in D.$$
Furthermore, equation (2.8) implies that
\be
\log \eta_3(\frac{t}{a})+ \log \eta_3(\frac{t}{b})=0, t\in D.
\ee
Let $\zeta(t)= \log \eta_3(t), t \in D.$ Equation (2.11) implies that
$$\zeta(\frac{t}{a})+\zeta(\frac{t}{b})=0, t \in D$$
or equivalently
$$\zeta(u)=-\zeta(\lambda u), u \in D^*$$
where $\lambda =\frac{a}{b}$ and $D^*= \{a^{-1}t, t \in D\}.$ Suppose that $\lambda =1.$ Then it follows that $\zeta(u)=0, u\in D^*$ which in turn implies that $\eta_3(t)=1, t\in D$ and hence $F_{Z_1}(t)=F_{S_1}(t), t \in D.$ If $\lambda \neq 1,$ applying Lemma 2 in Li and Zheng (2019), it follows that $\zeta(u)=0, u \in D^*$ under the additional condition of differentiability of the function $\zeta(u)$ and continuity of its derivative which implies that $F_{Z_1}(t)=F_{S_1}(t), t \in D.$ 
\vsp
This completes the proof of Theorem 2.2.
\vsp
We will now investigate the same problem when the constants $a>0, b<0,c>0$ and $d<0.$ 
\vsp
\NI{\bf Theorem 2.3:} Let $U$ and $V$ be random variables defined by the relations
$$U=\max(X,aZ_1,bZ_2);\;V=\max(Y,cZ_1,dZ_2)$$
where $X,Y,Z_1$ and $Z_2$ are independent random variables, $Z_1,Z_2$ are identically distributed and $a>0,b<0,c>0,d<0,$ are known  constants. Further suppose that the distribution functions of $X,Y,Z_1$ have the same support. Then, the joint distribution function of $(U,V)$ are connected by the distributions of $X,Y,Z_1$ and the distributions of $M,N,S_1$ by the equations through the equations (2.12) and (2.13) given below:
\be
F_X(t)=F_M(t) \frac{F_{S_1}(\frac{t}{a}) (1-F_{S_1}(\frac{t}{b}))}{F_{Z_1}(\frac{t}{a}) (1-F_{Z_1}(\frac{t}{b}))} 
\ee
and
\be
F_Y(t)=F_N(t) \frac{F_{S_1}(\frac{t}{c}) (1-F_{S_1}(\frac{t}{d}))}{F_{Z_1}(\frac{t}{c}) (1-F_{Z_1}(\frac{t}{d}))}.
\ee
\vsp
\NI{\bf Proof:} We follow the same notation as given in the proof of Theorem 2.2. It easy to see that that the joint distribution of The bivariate random vector $(U,V)$ is given by
\ben
G(t_1,t_2)&=&P(X\leq t_1,Y\leq t_2,Z_1\leq \frac{t_1}{a},Z_1\leq \frac{t_2}{c},Z_2\geq \frac{t_1}{b},Z_2\geq \frac{t_2}{d})\\\nonumber
& =& P(X\leq t_1,Y\leq t_2,Z_1\leq \min(\frac{t_1}{a},\frac{t_2}{c}),Z_2\geq \max(\frac{t_1}{b},\frac{t_2}{d}))\\\nonumber
&=& F_X(t_1)F_Y(t_2)F_{Z_1}(\min(\frac{t_1}{a},\frac{t_2}{c}))(1-F_{Z_1}(\max(\frac{t_1}{b},\frac{t_2}{d}))\\\nonumber
&=& F_M(t_1)F_N(t_2)F_{S_1}(\min(\frac{t_1}{a},\frac{t_2}{c}))(1-F_{S_1}(\max(\frac{t_1}{b},\frac{t_2}{d}))\\\nonumber 
\een
for all $t_1,t_2 \in D.$ Define $\eta_1(t), \eta_2(t)$ and $\eta_3(t)$ as defined in the proof of Theorem 2.2. Then it follows that
\be
\eta_1(t_1)\eta_2(t_2)\eta_3(\min(\frac{t_1}{a},\frac{t_2}{c}))\frac{ (1-F_{Z_1}(\max(\frac{t_1}{b},\frac{t_2}{d}))}{(1-F_{S_1}(\max(\frac{t_1}{b},\frac{t_2}{d}))}=1, t_1,t_2 \in D. 
\ee
Let $t_2 \raro \infty$ or to the right boundary point of the set $D.$ Then, it follows that,
\be
\eta_1(t_1) \eta_3(\frac{t_1}{a})\frac{(1-F_{Z_1}(\frac{t_1}{b}))}{(1-F_{S_1}(\frac{t_1}{b}))}=1 
\ee
observing that $\frac{t_2}{d} \raro -\ity$ as $t_2\raro \ity.$ Letting $t_1 \raro \ity,$ it follows that
\be
\eta_2(t_2)\eta_3(\frac{t_2}{c})\frac{(1-F_{Z_1}(\frac{t_2}{d}))}{(1-F_{S_1}(\frac{t_2}{d}))}=1.
\ee
These equations in turn show that
\be
F_X(t)=F_M(t) \frac{F_{S_1}(\frac{t}{a}) (1-F_{S_1}(\frac{t}{b}))}{F_{Z_1}(\frac{t}{a}) (1-F_{Z_1}(\frac{t}{b}))} 
\ee
and
\be
F_Y(t)=F_N(t) \frac{F_{S_1}(\frac{t}{c}) (1-F_{S_1}(\frac{t}{d}))}{F_{Z_1}(\frac{t}{c}) (1-F_{Z_1}(\frac{t}{d}))}.
\ee
\newsection{Identifiability by maxima for max-independent random variables}

{\bf Definition:} A finite collection of random variables $X_1,\dots,X_n$ is said to be {\it max-independent} if there exists a function $\beta(x_1,\dots,x_n)$ such that 
$$F(x_1,\dots,x_n)=F_1(x_1)\dots F_n(x_n)\beta(x_1,\dots, x_n), x_i\in R, 1\leq i \leq n$$
where $F(x_1,\dots,x_n)$ is the joint distribution of $(X_1,\dots,X_n)$ , $F_i(x)$ is the distribution function of $X_i$ for 
$1\leq i \leq n $ and $\beta(x_1,\dots,x_n)$ is a function taking values in the interval $(0,1]$ such that $\beta(x_1,\dots,x_n) \raro 1$ if $x_i\raro \ity$ for some $i, 1\leq i \leq n$ (cf. Prakasa Rao (2023)). The function $\beta(x_1,\dots,x_n)$ is called the {\it generator} of the random sequence $X_1,\dots,X_n.$ 
\vsp
Examples of sequence of max-independent random variables which are not independent are given in Prakasa Rao (2023). We will now generalize Theorem 2.2 to max-independent random variables. 
\vsp
\NI{\bf Theorem 3.1:} Let $X,Y, Z_1$ and $Z_2$ be max-independent random variables with generator $\beta (x_1,x_2,x_3,x_4).$ Define the random variables $U$ and $V$ by the relations
$$U=\max(X,aZ_1,bZ_2);\;V=\max(Y,cZ_1,dZ_2)$$
where $Z_1,Z_2$ are identically distributed and $a,b,c,d$ are known positive constants. Further suppose that the distribution functions of $X,Y,Z_1$ have the same support. Then, the joint distribution function of $(U,V)$ uniquely determines the distributions of $X,Y,Z_1$ if $a=b$ or if $a \neq b$ but the distribution function of $Z_1$ is differentiable and the derivative is continuous.
\vsp
\NI{\bf Proof:} Let $G(t_1,t_2)$ be the joint distribution function of $(U,V).$ Suppose that $(N,S_1,S_2)$ is another set of max-independent random variables with the generator $\beta (x_1,x_2,x_3,x_4)$ such that the joint distribution of the random vector
$$(\max(N,aS_1,bS_2), \max(M,cS_1,dS_2))$$
is the same as that of the random vector $(U,V).$ Let $F_X(.)$ denote the distribution function of the random variable $X.$ It can be seen that, for any$ -\ity<t_1,t_2<\ity,$
\bea
\;\;\;\;\\\nonumber
G(t_1,t_2)&=& P (\max(X,aZ_1,bZ_2)\leq t_1,\max((Y,cZ_1,dZ_2)\leq t_2)\\\nonumber
&=& P(X \leq t_1,aZ_1\leq t_1,bZ_2\leq t_2; Y\leq t_2,cZ_1\leq t_2, dZ_2\leq t_2)\\\nonumber
&=& P(X\leq t_1, Y\leq t_2, Z_1\leq \frac{t_1}{a}, Z_1\leq \frac{t_2}{c}, Z_2\leq \frac{t_1}{b}, Z_2\leq \frac{t_2}{d})\\\nonumber
&=& P(X\leq t_1, Y\leq t_2, Z_1\leq \min(\frac{t_1}{a},\frac{t_2}{c}), Z_2\leq \min(\frac{t_1}{b},\frac{t_2}{d}))\\\nonumber
&=& F_X(t_1)F_Y(t_2)F_{Z_1}(\min(\frac{t_1}{a},\frac{t_2}{c}))F_{Z_1}(\min(\frac{t_1}{b},\frac{t_2}{d}))\beta(t_1,t_2,\min(\frac{t_1}{a},\frac{t_2}{c}),\min(\frac{t_1}{b},\frac{t_2}{d})) \\\nonumber
\eea
since $X,Y Z_1,Z_2$ are max-independent random variables,$Z_1,Z_2$ are identically distributed with generator $\beta(x_1,x_2,x_3,x_4).$ By a similar argument it follows that
\bea
\;\;\;\;\\\nonumber
G(t_1,t_2)&=& P (\max(N,aS_1,bS_2)\leq t_1,\max((M,cS_1,dS_2)\leq t_2))\\\nonumber
&=& F_N(t_1)F_M(t_2)F_{S_1}(\min(\frac{t_1}{a},\frac{t_2}{c}))F_{S_1}(\min(\frac{t_1}{b},\frac{t_2}{d}))\beta(t_1,t_2,\min(\frac{t_1}{a},\frac{t_2}{c}),\min(\frac{t_1}{b},\frac{t_2}{d}) .\\\nonumber
\eea
Hence
\bea
F_X(t_1)F_Y(t_2)F_{Z_1}(\min(\frac{t_1}{a},\frac{t_2}{c}))F_{Z_1}(\min(\frac{t_1}{b},\frac{t_2}{d}))
\beta(t_1,t_2,\min(\frac{t_1}{a},\frac{t_2}{c}),\min(\frac{t_1}{b},\frac{t_2}{d}) \\\nonumber
=F_N(t_1)F_M(t_2)F_{S_1}(\min(\frac{t_1}{a},\frac{t_2}{c}))F_{S_1}(\min(\frac{t_1}{b},\frac{t_2}{d}))\beta(t_1,t_2,\min(\frac{t_1}{a},\frac{t_2}{c}),\min(\frac{t_1}{b},\frac{t_2}{d})  \\\nonumber 
\eea
for all $-\ity <t_1, t_2<\ity.$ Let $D\subset R$ be the common support of $X,Y,Z_1$ and $N,M,S_1.$ Then, for all $t\in D,$,
define
$$\eta_1(t)= \frac{F_X(t)}{F_M(t)}, \eta_2(t)= \frac{F_Y(t)}{F_N(t)},\mbox{and} \eta_3(t)= \frac{F_{Z_1}(t)}{F_{S_1}(t)}.$$
Equation (3.3) implies that
\be
\eta_1(t_1)\eta_2(t_2)\eta_3(\min(\frac{t_1}{a},\frac{t_2}{c}))\eta_3(\min(\frac{t_1}{b},\frac{t_2}{d}))=1,t_1, t_2 \in D. 
\ee
Let $t_1\raro \ity$ or $t_1$ tend to the right boundary point of the set $D$ in equation (3.4). From the properties of the distribution functions, it follows that 
$$\eta_2(t_2)\eta_3(\frac{t_2}{c})\eta_3(\frac{t_2}{d})=1, t_2\in D.$$
Hence
\be
\eta_2(t_2)=[\eta_3(\frac{t_2}{c})\eta_3(\frac{t_2}{d})]^{-1}, t_2\in D. 
\ee
Let $t_2 \raro \ity$ or $t_1$ tend to the right boundary point of the set $D$ in equation (3.4). From the properties of the distribution functions again, it follows that 
$$\eta_1(t_1)\eta_3(\frac{t_1}{a})\eta_3(\frac{t_1}{b})=1, t_1\in D$$ 
which implies that
\be
\eta_1(t_1)=[\eta_3(\frac{t_1}{a})\eta_3(\frac{t_1}{b})]^{-1}, t_2\in D. 
\ee
Combining the equations (3.4)-(3.6), it follows that
\be
\eta_3(\min(\frac{t_1}{a},\frac{t_2}{c}))\eta_3(\min(\frac{t_1}{b},\frac{t_2}{d}))= \eta_3(\frac{t_1}{a})\eta_3(\frac{t_1}{b})\eta_3(\frac{t_2}{c})\eta_3(\frac{t_2}{d}), t_1,t_2\in D.
\ee
Fix $t_1 \in D$ and let $t_2 \raro \ity$ or $t_2$ tend to the right boundary point of the set $D$. Then the expression on the left side of equation (3.6) tends to 1 and the expression on the right side tends to 
$$\eta_3(\frac{t_1}{a})\eta_3(\frac{t_1}{b}).$$
Hence
\be
\eta_3(\frac{t_1}{a})\eta_3(\frac{t_1}{b})=1, t_1\in D. 
\ee
This in turn implies that
\be
\eta_1(t_1)=1, t_1\in D. 
\ee
A similar analysis shows that 
\be
\eta_2(t_2)=1, t_2\in D 
\ee
by fixing $t_2\in D$ and letting $t_1\raro \ity$ or $t_1$ tend to the right boundary point of the set $D$.
In particular, it follows that
$$F_X(t)=F_M(t), t \in D$$
and
$$F_Y(t)=F_N(t), t \in D.$$
Furthermore, equation (3.8) implies that
\be
\log \eta_3(\frac{t}{a})+ \log \eta_3(\frac{t}{b})=1, t\in D.
\ee
Let $\zeta(t)= \log \eta_3(t), t \in D.$ Equation (3.11) implies that
$$\zeta(\frac{t}{a})+\zeta(t)(\frac{t}{b})=0, t \in D$$
or equivalently
$$\zeta(u)=-\zeta(\lambda u), u \in D^*$$
where $\lambda =\frac{a}{b}$ and $D^*= \{a^{-1}t, t \in D\}.$ Suppose that $\lambda =1.$ Then it follows that $\zeta(u)=0, u\in D^*$ which in turn implies that $\eta_3(t)=1, t\in D$ and hence $F_{Z_1}(t)=F_{S_1}(t), t \in D.$ If $\lambda \neq 1,$ applying Lemma 2 in Li and Zheng (2019), it follows that $\zeta(u)=0, u \in D^*$ under the additional condition of differentiability of the function $\zeta(u)$ and continuity of its derivative which implies that $F_{Z_1}(t)=F_{S_1}(t), t \in D.$ 
\vsp
This completes the proof of Theorem 3.1.
\vsp
\NI{\bf Acknowledgment:} This work was supported under the scheme ``INSA Honorary Scientist" at the CR Rao Advanced Institute of Mathematics, Statistics and Computer Science, Hyderabad 500046, India.  
\vsp
\NI{\bf References}
\begin{description}

\item Klebanov, L. 1973. Reconstructing the distributions of the components of a random vector from distributions of certain statistics, {\it Mathematical Notes}, 13:71-72.

\item Kotlarski, I. 1978. On some characterization in probability by using minima and maxima of random variables, {\it Aequationes Mathematicae}, 17:77-82.

\item Li, Siran., and Zheng, Xunjie. 2019. A generalization of Lemma 1 in Kotlarski (1967). Preprint.

\item Prakasa Rao, B.L.S. 1992. {\it Identifiability in Stochastic Models: Characterization of Probability Distributions}. New York: Academic Press.

\item Prakasa Rao , B.L.S. 2023. On some characterizations of probability distributions based on maxima or minima of some families of dependent random variables, {\it Jour. Indian Soc. Probab. and  Statist.}, https://doi.org/10.10007/s41096-024-00180-1.

\end{description}

\end{document}